\documentclass[runningheads]{llncs}
\usepackage{graphicx}
\usepackage{amssymb}
\def\cro{{\mbox {\sc cr}}}
\newtheorem{observation}{Observation}

\begin{document}

\title{Improvement on the Crossing Number of Crossing-Critical Graphs\thanks{Supported
    by National Research, Development and Innovation Office, NKFIH,
K-131529 and the Higher Educational Institutional Excellence Program 2019, the grant of the 
Hungarian Ministry for Innovation and Technology (Grant Number: NKFIH-1158-6/2019).}}
%
%
\author{J\'anos Bar\'at\inst{1}\orcidID{0000-0002-8474-487X} \and
G\'eza T\'oth\inst{2,3}\orcidID{0000-0003-1751-6911}}  

\authorrunning{J. Bar\'at and G. T\'oth}
%
\institute{Department of Mathematics, University of Pannonia, Veszpr\'em \email{barat@mik.pannon.hu} \and
 Alfr\'ed R\'enyi Institute of Mathematics, \email{geza@renyi.hu} \and  
Budapest University of Technology and Economics, SZIT}

\maketitle              
\begin{abstract}
The crossing number of a graph $G$ is
the minimum number of edge crossings over all drawings of $G$
in the plane.
A graph $G$ is $k$-crossing-critical if
its crossing number is at least $k$, but
if we remove any edge of $G$, its crossing number drops below $k$.
There are examples of $k$-crossing-critical graphs that do not have drawings with exactly $k$ crossings.
Richter and Thomassen proved in 1993 that if $G$ is 
$k$-crossing-critical, then its crossing number is at most $2.5k+16$.
We improve this bound to $2k+6\sqrt{k}+47$.
\end{abstract}

\section{Introduction}

The crossing number $\cro(G)$ of a graph $G$ is the
minimum number of edge crossings over all drawings of $G$
in the plane.
In the optimal drawing of $G$, crossings are not necessarily
distributed uniformly on the edges. 
Some edges can be more ``responsible'' for the crossing number than others.
For any positive integer $k$, there exists a graph $G$ whose crossing number is $k$,
but it has  an edge $e$ such that $G-e$ is planar.

On the other hand, Richter and Thomassen \cite{RT93} (Section 3.) conjectured that if $\cro(G)=k$, then
$G$ contains an edge $e$ such that  $\cro(G-e)\ge k-c\sqrt{k}$ for some constant $c$.
They observed that this bound would be optimal, as shown, e.g., by the graph $K_{3,n}$.
They managed to prove a much weaker bound, namely, 
if $\cro(G)=k$, then
$G$ contains an edge $e$ such that  $\cro(G-e)\ge 2k/5-8$.

A graph $G$ is {\em $k$-crossing-critical}
if $\cro(G)\ge k$, but $\cro(G-e)<k$ for any edge $e$ of $G$.

The structure and properties of
crossing-critical graphs are fundamental in the study of crossing numbers. 
It is easy to describe $1$-crossing-critical graphs, and there is an almost
complete description of $2$-crossing-critical graphs \cite{BORS16}.
For $k>2$, a description of $k$-crossing-critical
graphs seems hopeless at the moment.

It has been proved recently, that the bounded maximum degree conjecture for 
$k$-crossing-critical graphs holds for $k\le 12$ and does not hold for
$k>12$ \cite{BDH19}. More precisely, there is a constant $D$ with the property
that for every $k\le 12$, every $k$-crossing-critical graph has maxmimum degree
at most $D$, and for every $k>12$, $d\ge 1$, there is a
$k$-crossing-critical graph with maximum degree at least $d$.

\smallskip

We rephrase the result and conjecture of
Richter and Thomassen \cite{RT93} as follows.
They conjectured that if $G$ is $k$-crossing-critical, then $\cro(G)\le k+c'\sqrt{k}$ for some $c'>0$
and this bound would be optimal.
They proved that if $G$ is $k$-crossing-critical, then $\cro(G)\le 2.5k+16$.
This result  has been improved
in two special cases.

Lomel\'{\i} and Salazar \cite{LS06} proved that for any $k$ there is an $n(k)$ such that
if $G$ is $k$-crossing-critical and has at least
$n(k)$ vertices, then $\cro(G)\le 2k+23$.

Salazar \cite{S00} proved that if $G$ is $k$-crossing-critical and
all vertices of $G$ have degree at least $4$, then $\cro(G)\le 2k+35$.

It is an easy consequence of the Crossing Lemma \cite{A19} that
if the average degree in a $k$-crossing-critical graph is large, then
its crossing number is close to $k$ \cite{FT08}. More precisely,
if $G$ is $k$-crossing critical and it has at least $cn$ edges, where 
$c\ge 7$, 
then
$\cro(G)\le kc^2/(c^2-29)$.


In this note, we obtain a general improvement.

%
\begin{theorem}\label{main}
  For any $k>0$, if $G$ is a $k$-crossing-critical multigraph, 
  then $\cro(G)\le 2k+6\sqrt{k}+47$.
 \end{theorem}

\medskip

We need a few definitions and introduce now several
parameters for the proof. We also list them at the end of the paper.

Let $G$ be a graph. 
We call a pair $(C, v)$, where $C$ is a cycle of $G$ and $v$ is a vertex of $C$,
the {\em cycle $C$ with special vertex $v$}. 
(The special vertex meant to be a vertex with large degree.)
When it is clear from the context, which one is the special vertex,
we just write $C$ instead of $(C, v)$.

Suppose $C$ is a cycle with special vertex $v$.
Let $x$ be a vertex of $C$. 
An edge, adjacent to $x$ but not in $C$, 
is {\em hanging from $x$} in short.
Let $l(C)=l(C, v)$ be the {\em length of $C$}, 
that is, the number of its edges.
For any vertex $x$, let $d(x)$ denote the degree of $x$. 
Let $h(C)=h(C, v)=\sum_{u\in C, u\neq v}(d(u)-2)$, 
that is,
the total number of hanging edges
from all non-special vertices of $C$ (with multiplicity).

A set of edges is {\em independent} if no two of them have a common endvertex.

\section{The Proof of Richter and Thomassen}

In \cite{RT93}, the most important tool in the proof was the following technical result.
In this section we review and analyze its proof. The algorithmic argument finds a cycle $C$ recursively  
such that  $h(C)$ is small.

\medskip

\noindent {\bf Theorem 0.} \cite{RT93}
{\em Let $H$ be a simple graph with minimum degree at least $3$.
Assume that $H$ has a set $E$ of $t$ edges such that $H-E$ is planar. Then
$H$ has a cycle $K$ with special vertex $v$ such that $h(K)\le t+36$.}

\medskip

\noindent {\bf Proof of Lemma~0.}
The proof is by induction on $t$. 
The induction step can be  considered
as a process, which constructs a graph $H^*$ from
graph $H$, and
cycle $K$ of $H$, either directly, or from cycle $K^*$ in $H^*$.
For convenience, for any planar graph $H$, define $H^*=\emptyset$.
In the rest of the paper we refer to this as the Richter-Thomassen procedure.
The statement of Theorem 0 for $t=0$ is the following.

\smallskip

\noindent {\bf Lemma 0.} \cite{RT93}
{\em Let $H$ be a simple planar graph with minimum degree at least $3$.
Then
$H$ has a cycle $K$ with special vertex $v$ such that $l(K)\le 5$ and $h(K)\le 36$.}

\smallskip

Here we omit the proof of Lemma 0.
Suppose now that $t>0$ and we have already shown Theorem 0 for smaller values of $t$.
Let $H$ be a simple graph with minimum degree at least $3$.
Assume that $H$ has a set $E$ of $t$ edges such that $H-E$ is planar and let $e=uw\in E$.
Let $H'=H-e$. We distinguish several cases.

\bigskip

\begin{itemize}

  \item[\underline{1.}]
    $H'$ has no vertex of degree $2$.
By the induction hypothesis, $H'$ has a cycle $K^*$ with a special vertex $v$ such that
$h(K^*)\le t+35$. If $e$ is not a chord of $K^*$, then $K=K^*$ with the same special vertex satisfies the conditions for
$H$. Let $H^*=H'$. 

If $e$ is a chord of $K^*$, then $K^*+e$ determines two cycles,
and it is easy to see that either one satisfies the conditions. So, let $K$ be one of them. 
If $K$, contains $v$, then $v$ remains the special vertex.
If $K$ does not contain $v$, then we can choose the special vertex of $K$ arbitrarily.
Let $H^*=H'$.

\end{itemize}

\smallskip

\begin{itemize}

\item[\underline{2.}]
  $H'$ has a vertex of degree $2$.
Clearly, only $u$ and $w$ can have degree $2$. 
Suppress vertices of degree $2$. That is, for each vertex of degree 2,  
remove the vertex and connect its neighbors by an edge.
Let $H''$ be the
resulting graph. It can have at most two sets of parallel edges.

\medskip

\begin{itemize}


\item[\underline{2.1.}] 
$H''$ has no parallel edges. By the induction hypothesis, $H''$ 
contains a cycle $K^*$ with a special vertex $v$
such that $h(K^*)\le t+35$.
It corresponds to a cycle $K'$ in $H$. 
Let $H^*=H''$.

\begin{itemize}

\item[\underline{2.1.1.}]
  The edge $e$ is not incident with $K'$. In this case, $K=K'$
satisfies the conditions, with the same special vertex as $H^*$. 

\item[\underline{2.1.2.}]
  The edge $e$ has exactly one endvertex on $K'$. 
In this case, let $K=K'$
with the same special vertex. 
Now $h(K)=h(K^*)+1\le t+36$ and we are done.

\item[\underline{2.1.3.}]
  The edge $e$ has both endvertices on $K'$.
Now, just like in Case 1, $K'+e$ determines two cycles and
it is easy to see that either one satisfies the conditions.
If the new cycle contains $v$, then it will remain the special vertex, if
not, then we can choose the special vertex arbitrarily.

\end{itemize}

\medskip


\item[\underline{2.2.}] 
  $H''$ has one set of parallel edges. 
Let $x$ and $y$ be the endvertices of the parallel edges.
We can assume that one of the $xy$ edges in $H''$ 
corresponds to the path $xuy$ in $H$ and $H'$.
Clearly, $d(u)=3$. 

\begin{itemize}

\item[\underline{2.2.1.}]
Another  $xy$ edge in $H''$ 
corresponds to the path $xwy$ in $H$ and $H'$.
In this case  $d(u)=d(w)=3$, so for the cycle 
$K=uxw$ with special vertex $x$ we have $h(K)\le 2$ and we are done.
Let $H^*=\emptyset$.
We do not define $K^*$ in this case.

\item[\underline{2.2.2.}]
No $xy$ edge in $H''$ 
corresponds to the path $xwy$ in $H$ and $H'$ and 
either $d(x)\le 37+t$
  or $d(y)\le 37+t$.
Assume that $d(x)\le 37+t$, the other case is treated analogously.
Since there were at least two $xy$ edges in $H''$, $H$ contains the edge $xy$. 
For the cycle $K=uxy$, with special vertex $y$, we have $h(K)\le 35+t+1$,
so we are done.
Let $H^*=\emptyset$.


\item[\underline{2.2.3.}]
No $xy$ edge in $H''$ 
corresponds to the path $xwy$ in $H$ and $H'$ and 
both 
$d(x), d(y) > 37+t$.
Replace
the parallel edges by a single $xy$ edge in $H''$. In the resulting graph $H^*$,
we can apply the induction hypothesis and get a cycle $K^*$ with special
vertex $v$ such that
$h(K^*)\le 35+t$. 
Now $K^*$ cannot contain both $x$ and $y$ and if it
contains either one, then it has to be the special vertex.
Therefore, the cycle $K$ in $H$, correspondig to $K^*$, 
with the same special vertex, satisfies the conditions, since the only edge that can increase $h(K)$ is $e$,
and $e$ is not a chord of $K$.

\end{itemize}


\medskip

\item[\underline{2.3.}] 
  $H''$ has two sets of parallel edges, $xy$ and $ab$ say. 
Now $H$ contains the edges $xy$ and $ab$. 
We can assume by symmetry that $H$ contains the paths  $xuy$ and  $awb$.
Also $d(u)=d(w)=3$ in $H$.

\begin{itemize}

\item[\underline{2.3.1.}]
  At least one of $a$, $b$, $x$, $y$ has degree at most $37+t$ in $H$.
Assume that $d(x)\le 37+t$, the other cases are  treated analogously.
For the cycle $K=uxy$ with special vertex $y$, we have $h(K)\le 35+t+1$,
so we are done.
Let $H^*=\emptyset$.

\item[\underline{2.3.2.}]
  $d(x)$, $d(y)$, $d(a)$, $d(b)>37+t$. 
Replace
the parallel edges by single edges $xy$ and $ab$ in $H''$. In the resulting
graph $H^*$,
we can apply the induction hypothesis and get a cycle $K^*$ with special
vertex $v$ such that
$h(K^*)\le 35+t$. However, $K^*$ can contain at most one of $x$, $y$, $a$,
and $b$, and if it contains one, that
has to be the special vertex.
Therefore, the cycle $K$ in $H$, correspondig to $K^*$
with the same special vertex satisfies the conditions.

\end{itemize}
\end{itemize}
\end{itemize}

This finishes the proof of Theorem 0. 
$\Box$





\section{Proof of Theorem~\ref{main}.}

The main idea in the proof of Richter and Thomassen \cite{RT93} is the following.
Suppose that $G$ is $k$-crossing-critical. Then it has at most $k$ edges whose removal
makes $G$ planar. Then by Theorem 0, we find a cycle $C$ with special vertex $v$ such that $h(C)\le k+36$.
Let $e$ be an edge of $C$, adjacent to $v$. We can draw $G-e$ with at most $k-1$ crossings.
Now we add the edge $e$, along $C-e$, on the ``better'' side.
We get additional crossings from the crossings on $C-e$, and from the hanging edges and we can bound both.

Our contribution is the following.
Take a ``minimal'' set of edges, whose removal makes $G$ planar. Clearly, this set would contain at most $k$ edges.
However, we have to define ``minimal'' in a slightly more complicated way, 
but still, our set contains at most $k+\sqrt{k}$ edges.
We carefully analyze the proof of Richter and Thomassen, extend it with some operations, and find a 
cycle $C$ with special vertex $v$ such that (roughly) $l(C)+h(C)/2\le k+6\sqrt{k}$.
Now, do the redrawing step. If $h(C)$, or the number of crossings on $C-e$ is small,
then we get an improvement immediately.
If both of them are large, then  $l(C)$ is much smaller than the number of crossings on $C-e$.
But in this case,  we can remove the edges of $C$, and get rid of many crossings. This way, we can
get a bound on the ``minimal'' set of edges whose removal
makes $G$ planar.

\smallskip

As we will see, for the proof we can assume that $G$ is simple and all vertices have degree at least $3$.
But if we want to prove a better bound, say, $\cro(G)\le (2-\varepsilon)k+o(k)$,
then we cannot prove that the result for simple graphs implies the result for multigraphs. Therefore,
the whole proof collapses. Moreover, even if we could assume without loss of generality
that $G$ is simple, we still cannot go below the constant $2$ with our method. We cannot rule out the possibility
that all (or most of the) $k-1$ crossings are on $C-e$.

\bigskip

\noindent {\bf Proof of Theorem \ref{main}.}
Suppose that $G$ is $k$
-crossing-critical.
Just like in the paper of Richter and Thomassen \cite{RT93}, we can assume that
$G$ is simple and all vertices have degree at least $3$.
We sketch the argument.

If $G$ has an isolated vertex, we can remove it from $G$.
Suppose that a vertex $v$ of $G$ has degree $1$. Then $\cro(G)=\cro(G-v)$,
contradicting crossing criticality.
Suppose now that $v$ has degree $2$.
We can suppress $v$ (remove it and connect its neighbors by an edge).
The resulting (multi)graph is still $k$-crossing-critical and
has the same crossing number as $G$.

Clearly, $G$ cannot contain loops, as adding or removing a loop does not
change the crossing number.
Finally, suppose that $e$ and $f$ are parallel edges, both connecting $x$
and $y$.
Since  $G$ is $k$-crossing-critical, we have $\cro(G-e)\le k-1$. Take a
drawing of the graph $G-e$ with at most $k-1$ crossings.
Add the edge $e$, drawn very close to $f$. The obtained drawing of $G$ has
at most $2k-2$ crossings.

So, we assume in the sequel, that $G$ is simple and all vertices have degree 
at least $3$.

\bigskip

Let $k'$ 
be the smallest integer with the property that we can remove $k'$ edges from
$G$ so that the remaining graph is planar. 
Define the function $f(x, y)=\sqrt{k}x+y$.

Let $(t,t')$ 
be the pair of numbers that minimizes the function
$f(t, t')=\sqrt{k}t+t'$ subject to the following property:
There exists a set $E$ of $t$ edges such that $G-E$ is planar, and the set 
$E$ contains at most $t'$ independent edges. 
In the next lemma, part (i) is from \cite{RT93}, we repeat it here for completeness.

\begin{lemma} \label{lem1}
The following two statements hold.\\
{\rm (i) \cite{RT93}} $k'\le k$, and\\ 
{\rm (ii)} $t\le k'+\sqrt{k}$.
\end{lemma}

\noindent {\bf Proof of Lemma~\ref{lem1}.} (i) Since $G$ is $k$-crossing-critical,  $G-e$
can be drawn with at most $k-1$ crossings for any edge $e$. 
Remove one of the edges from each crossing in such a drawing.
We removed at most $k$ edges in total and got a planar graph.
(ii) Let $E'$ be a set of $k'$ edges such that $G-E'$ is planar. 
Suppose that $E'$ contains at most $k''$ independent edges. Now $k''\le k'$.
By the choice of $(t,t')$,
$f(t,t')\le f(k',k'')$. Consequently,
$\sqrt{k}t\le \sqrt{k}t+t'=f(t,t')\le f(k',k'')\le \sqrt{k}k'+k'$.
Therefore,   $t\le k'+k'/\sqrt{k}\le k'+\sqrt{k}$. 
$\Box$

\medskip

Now set $E=\{ e_1, e_2, \ldots, e_t \}$, where $E$ contains at most 
$t'$ independent edges, and $G-E$ is planar.
Apply the Richter-Thomassen procedure recursively starting with $H_0=G$.
We obtain a sequence of graphs $H_0, H_1, \ldots, H_s$, ($s\le t$)
such that 
for $0\le i\le s-1$, $H_i^*=H_{i+1}$, and $H_s^*=\emptyset$. 
The procedure stops with graph $H_s$, where we obtain a cycle $C_s$ either directly,
in cases 2.2.1, 2.2.2, and 2.3.1, or by Lemma 0, when $H_s$ is planar. In all cases, $l(C_s)\le 5$. 
Following the procedure again, we also obtain 
cycles $C_{s-1}, \ldots, C_0$ of
 $H_{s-1}, \ldots, H_0$ respectively
such that
 $0\le i\le s-1$, $C_i^*=C_{i+1}$. 
Let $C_0=C$ with special vertex $v$.

\begin{lemma}\label{lem2}
 There is a cycle $K$ of $G$ such that
$l(K)+h(K)/2\le t+5\sqrt{k}+48$.
\end{lemma}

\noindent {\bf Proof of Lemma~\ref{lem2}.} The cycle $K$ will be either $C$, or a slightly 
modified version of $C$.
It is clear from the procedure that $C$ does not have a chord in $G$ since
we always choose $C$ as a minimal cycle. 

Consider the moment of the procedure, when we get cycle $K$ from $K^*$.
All hanging edges of $K$ correspond to a hanging edge of $K^*$, with the possible exception of $e=uw$.
Therefore, if we get a new hanging edge $e$,
then $e\in E$.
Taking into account the initial cases in the procedure, that is,
when we apply Lemma 0, or we have Cases 2.2.1, 2.2.2, or 2.3.1,
we get the following easy observations.
We omit the proofs. 

\begin{observation} \label{konnyu1}
{\rm (i)} All but at most $36$ edges of $G-C$
adjacent to a non-special vertex of $C$ are in $E$.\\
{\rm (ii)} For all but at most $4$ non-special vertices $z'$ of $C$, 
all edges of $G-C$ incident to $z'$, are in $E$. $\Box$
\end{observation}

\begin{figure}[ht]
\begin{center}
\scalebox{0.45}{\includegraphics{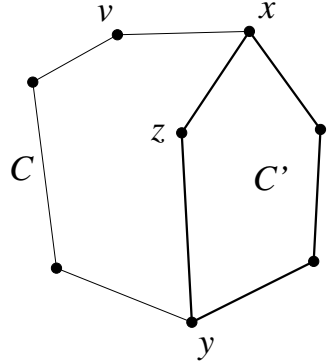}}
\caption{Cycles $C$ and $C'$.}\label{atloskor}
\end{center}
\end{figure}

Suppose that $l(C)>t'+6$.
Consider $t'+5$ consecutive vertices on $C$, none of them being the special vertex $v$.
By Observation~\ref{konnyu1}~(ii), for at least $t'+1$ of them, all hanging edges are in $E$. 
Consider one of these hanging edges at each of these $t'+1$ vertices. 
By the definition of $t'$, these
$t'+1$ edges cannot be independent, at least two of them have a common endvertex, which is not on $C$.
Suppose that $x, y\in C$, $z\not\in C$, $xz, yz\in E$.
Let $a$ be the $xy$ arc of $C$, which does not contain the special vertex $v$. Take two consecutive neighbors of $z$. Assume for simplicity,
that they are $x$ and $y$. 
Let the cycle $(C',z)$ be formed by arc $a$ of $C$, together with the path $xzy$. 
See {\bf Fig.~\ref{atloskor}}. The cycle $C'$ does not have a chord in $G$.
We have $l(C')\le t'+6$, $h(C')\le h(C)$.
The edges $zx$ and $zy$ are the only new hanging edges of $C'$ from a non-special vertex. They might not be in $E$, therefore,
the statement of Observation~\ref{konnyu1} holds in a slightly weaker form.

\begin{observation} \label{konnyu2}
{\rm (i)} All but at most $38$ edges of $G-C'$
adjacent to a non-special vertex of $C'$ are in $E$.\\
{\rm (ii)} For all but at most $6$ non-special vertices $z'$ of $C'$, 
all edges of $G-C'$ incident to $z'$, are in $E$. $\Box$
\end{observation}

\smallskip

Let cycle $K=C$, if $l(C)\le t'+6$, and let $K=C'$, if $l(C)> t'+6$. In both cases, for the rest of the proof, let $v$ denote the special vertex of $K$.
Let $h=h(K)$, $l=l(K)$.
We have 
\begin{equation} \label{eqlt}
 l\le t'+6.
\end{equation}

The cycle $K$
does not have a chord. 
In particular, 
none of 
$e_1, e_2, \ldots, e_t$ can be a chord of $K$.
Now we partition $E$ into three sets,
$E=E_p\cup E_q\cup E_m$, where 
$E_p$ 
is the 
subset of edges of $E$, 
which have exactly one endvertex on $K$ (these are the hanging edges in $E$),
$E_q=E\cap K$, 
$E_m$ is the subset of edges of $E$, which do not have an endvertex on $K$.
Let $p=|E_p|$, $q=|E_q|$, $m=|E_m|$. 
Let $p'$ denote the number of edges of $E_p$ hanging from the special vertex $v$.
By definition, 
\begin{equation} \label{eq1}
 t=p+q+m
\end{equation}

and 
$p\ge p'$.

It follows from Observations
\ref{konnyu1} (i) and \ref{konnyu2} (i) that 
\begin{equation} \label{eqhp}
 h\ge p-p'\ge h-38.
\end{equation}

Therefore,
$$h+q+m\le p+q+m+38=t+38.$$

Since all vertices have degree at least $3$, and $K$ does not have a chord,

$$h\ge l-1.$$

\medskip

\begin{figure}[ht]
\begin{center}
\scalebox{0.45}{\includegraphics{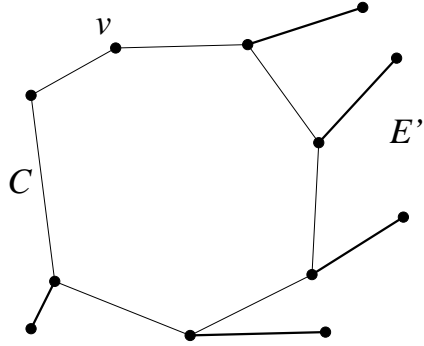}}
\caption{The edge set $E'$.}\label{tuskeskor}
\end{center}
\end{figure}

Now, 
at each vertex $x$ of $C$, where {\em all} hanging edges belong to $E_p$, take one such edge.
The set of these edges is $E'$.
By Observations
\ref{konnyu1} (ii) and \ref{konnyu2} (ii), $|E'|\ge l-7$. 
See {\bf Fig.~\ref{tuskeskor}}.
Let $F=E_p\cup E(K)\cup E_m- E'$ where $E(K)$ is the set of edges of $K$.
Since $F\cup E'\supseteq E$, $G'=G- (F\cup E')$ is a planar graph. Let $G''=G'\cup E'=G- F$.
In $G''$, each edge of $E'$ has an endvertex of degree one. 
Therefore, we can add all edges of $E'$ to $G'$ without losing planarity.
Consequently, the graph $G''=G'\cup E'=G- F$ is planar.

The set $F$ has at most 
$|F|\le p+l+m-(l-7)=p+7+m\le t+7$ edges by (\ref{eq1}).
That is
\begin{equation}\label{eqF<t+5}
 |F|\le t+7
\end{equation}

Let $F'\subseteq F$ be a maximal set of independent edges in $F$.
To estimate $|F'|$,
observe that apart from the edges in $E_m\cap F'$,
all edges in $F'$ are adjacent to a vertex of $C$.
Moreover, at most $p-p'-(l-7)+1$ of them have one vertex on $K$, the rest have two. 
Therefore, 
$$|F'|\le p-p'-(l-7)+1+(l-(p-p'-l+8))/2+m=(p-p')/2+4+m.$$
By the choice of the pair $(t,t')$,
$\sqrt{k}t+t'\le \sqrt{k}|F|+|F'|$.\\
Therefore,
$(p-p')/2+4+m\ge |F'|\ge \sqrt{k}t+t'-\sqrt{k}|F|\ge t'-5\sqrt{k}$, using (\ref{eqF<t+5}).

\bigskip

Now evoking (\ref{eqhp}):
$h/2+4+m\ge (p-p')/2+4+m\ge t'-5\sqrt{k}$. 

Therefore, 
$h/2+m\ge t'-4-5\sqrt{k}\ge l-10-5\sqrt{k}$ by (\ref{eqlt}).

Summarizing, we have 

$$h+m\le t+38, 
\ \  h/2+m\ge l-10-5\sqrt{k},$$
which implies 

$$m\le t-h+38, \ \  h/2+t-h+38\ge l-10-5\sqrt{k},$$
and finally
$$t+5\sqrt{k}+48\ge l+h/2.$$

This concludes the proof of Lemma~\ref{lem2}. $\Box$

\medskip

Now we can finish the proof of Theorem~\ref{main}. 
By Lemma~\ref{lem2}, we have a cycle $K$
in $G$ with special vertex $v$ such that
$h(K)/2+l(K)\le t+5\sqrt{k}+48$.
Let $e$ be an edge of $K$ adjacent to $v$. 
Since $G$ was $k$-crossing-critical,
the graph $G- e$
can be drawn with at most $k-1$ crossings. 
Let us consider such a drawing $D$. Let $h=h(K)$, $l=l(K)$.

Suppose the path $K- e$ has $cr$ crossings in $D$.
Remove the edges of $K$
from the drawing, and one edge from each crossing not on $K- e$.
Together with $e$, we removed at most
$k+l-cr$  edges from $G$ to get a planar graph. 
Therefore, 
$l+k-cr\ge k'.$ Combining it with Lemma~\ref{lem1} (ii) we have
$$l+k-cr\ge t-\sqrt{k}.$$
Consequently, 
$l+k-cr\ge t-\sqrt{k}\ge l+h/2-48-6\sqrt{k}$ by Lemma~\ref{lem2}. That is
$$k+6\sqrt{k}+48\ge cr+h/2.$$

\begin{figure}[h]
\begin{center}
\scalebox{0.45}{\includegraphics{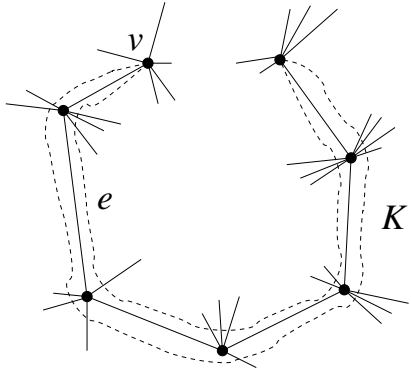}}
\caption{Adding the missing edge $e$.}\label{drawe}
\end{center}
\end{figure}

Consider the drawing $D$ of $G-e$.
We can add the missing edge $e$
drawn along the path $K- e$ on either side. See {\bf Fig.~\ref{drawe}}.
The two possibilities {\em together} create at most $h+2cr$ crossings. 
Choose the one which creates fewer crossings.
That makes at most $h/2+cr$ crossings.

Since
$k+6\sqrt{k}+48\ge cr+h/2$, 
 we can add $e$ with at most
$k+6\sqrt{k}+48$
additional crossings.
Hence $\cro(G)\le 2k+6\sqrt{k}+47$. $\Box$

\section*{Notations}

Here we give a list of the parameters and their definitions, used in the proof.

\bigskip

$k$: $G$ is $k$-crossing-critical

$k'$: the smallest integer with the property that we can remove $k'$ edges from
$G$ so that the remaining graph is planar.

$(t,t')$: the pair of numbers that minimizes the function $f(t, t')=\sqrt{k}t+t'$ subject to the following property:
There exists a set $E$ of $t$ edges such that $G-E$ is planar, and the set 
$E$ contains at most $t'$ independent edges. 

$p=|E_p|$: the number of edges in $E$ that have exactly one endvertex on $C$.

$q=|E_q|$: the number of edges in $E\cap C$.

$m=|E_m|$: the number  of edges in $E$ that do not have an endvertex on $C$.

$p'$: the number of edges of $E_p$ hanging from the special vertex $v$ of $C$.

$h=h(C)=h(C, v)=\sum_{u\in C, u\neq v}(d(u)-2)$, the total number of hanging edges
from all non-special vertices of $C$ (with multiplicity).

$l=l(C)$: the length of $C$.

\end{document}